\documentclass[11pt, reqno]{amsart}
\usepackage{amsmath, amssymb, amsthm, graphicx,marvosym}
\usepackage{cancel}
\usepackage[nobysame]{amsrefs}[11pts]

\newtheorem{theorem}{Theorem}[section]
\newtheorem{proposition}[theorem]{Proposition}
\newtheorem{lemma}[theorem]{Lemma}

\newtheorem{remark}[theorem]{Remark}

\makeatletter \@addtoreset{equation}{section} \makeatother

\newcommand{\pp}{\partial_+}
\newcommand{\pn}{\partial_-}
\def\tilde{\widetilde}

\newcommand{\ve}{\varepsilon}

\newcommand{\beq}{\begin{equation}}
\newcommand{\eeq}{\end{equation}}

\def\TS{\textstyle}
\def\com#1{\quad{\textrm{#1}}\quad}
\def\eq#1{(\ref{#1})}
\def\nn{\nonumber}
\def\ol{\overline}

\def\dbyd#1{\TS{\frac{\partial}{\partial#1}}}

\begin{document}

\title[Optimal bound on density for compressible Euler equations]{Optimal time-dependent lower bound on density for classical solutions of 1-D compressible Euler equations}

\author{Geng Chen}
\address{School of Mathematics,
Georgia Institute of Technology, Atlanta, GA 30332 USA ({\tt gchen73@math.gatech.edu}).}

\date{\today}

\begin{abstract}
{\small
For the compressible Euler equations, even when initial data are uniformly away from vacuum, solutions can approach vacuum in infinite time.   Achieving sharp lower bounds of density is crucial in the study of Euler equations.
In this paper, for the initial value problems of isentropic and full Euler equations in one space dimension, 
assuming the initial density has positive lower bound, we prove that density functions in classical solutions 
have positive lower bounds in the order of $\textstyle O(1+t)^{-1}$ and  $\textstyle O(1+t)^{-1-\delta}$ for any $\textstyle 0<\delta\ll 1$, respectively, where $t$ is time. The orders of these bounds are optimal or almost optimal, respectively.  Furthermore, for classical solutions in Eulerian coordinates $(y,t)\in\mathbb{R}\times[0,T)$, we show velocity $u$ satisfies that $u_{y}(y,t)$ is uniformly bounded from above by a constant independent of $T$, although $u_{y}(y,t)$ tends to negative infinity when gradient blowup happens, i.e. when shock forms, in finite time.
}

\end{abstract}
\maketitle

2010\textit{\ Mathematical Subject Classification:} 76N15, 35L65, 35L67.


\section{{\bf Introduction}}
The compressible Euler equations in Lagrangian coordinates in one space dimension are
\begin{align}
\tau_t-u_x&=0\,,\label{lagrangian1}\\
u_t+p_x&=0\,,\label{lagrangian2}\\
\textstyle\big(\frac{1}{2}u^2+e\big)_t+(u\,p)_x&=0\,, \label{lagrangian3}
\end{align}
where $\rho$ is the density, $\tau=\rho^{-1}$ is the specific volume,
$p$ is the pressure, $u$ is the velocity, $e$ is the  specific  internal energy, $t\in\mathbb{R}^+$ is the time
and $x\in\mathbb{R}$ is the spatial coordinate.  
The compressible Euler equations are widely used, especially in the gas dynamics. 
The classical solutions for the compressible Euler equations in Lagrangian and Eulerian coordinates are equivalent \cite{Dafermos2010}. 

For simplicity, in this paper, we only consider the case when the gas is ideal polytropic, in which
\beq
   p=K\,e^{\frac{S}{c_v}}\,\tau^{-\gamma} \com{with adiabatic gas constant} \gamma>1\,,
\label{introduction 3}
\eeq
and 
\[
e=\frac{p\tau}{\gamma-1}\,,
\]
where $S$ is the entropy,
$K$ and $c_v$ are positive constants, c.f. \cite{courant} or \cite{smoller}.  
For $C^1$ solutions,
it follows that (\ref{lagrangian3}) is equivalent to the conservation of entropy \cite{smoller}:
\beq
   S_t=0\,,
\label{s con}
\eeq
hence 
\[S(x,t)\equiv S(x,0)\doteq S(x).\]

If the entropy is constant, the flow is isentropic, then
(\ref{lagrangian1}) and (\ref{lagrangian2}) become a closed system,
known as the $p$-system:
\begin{align}
\tau_t-u_x&=0\,,\label{p1}\\
u_t+p_x&=0\,,\label{p2}
\end{align}
with
\beq\label{p3}
   p=K\,\tau^{-\gamma}\,,\qquad \gamma>1,
\eeq
where, without loss of generality, we still use $K$ to denote the constant in pressure.

In this paper, we consider the classical solutions of initial value problems for full Euler equations
(\ref{lagrangian1}), (\ref{lagrangian2}), (\ref{introduction 3}) and (\ref{s con}) with initial data
$\big(u(x,0), \tau(x,0), S(x,0)\big)$ and isentropic Euler equations (\ref{p1})$\sim$(\ref{p3}) with 
initial data $\big(u(x,0),  \tau(x,0)\big)$. 
We consider the large data problem, which means that there is no restriction on the size of the solutions.

Toward a large data global existence of BV solutions
for the compressible Euler equations, which is a major open problem in the field of
hyperbolic conservation laws, one of the main challenges is the possible
degeneracy when density approaches zero. In fact, a solution loses its strict hyperbolicity as density approaches zero. See \cites{BCZ,CJ,ls} for analysis and examples showing these difficulties.
Therefore, the sharp information on the time decay of density lower bound is critical in the study of compressible Euler equations. Furthermore, the time-dependent lower bound on density for classical solutions can be used to study the shock formation and life-span of classical solutions.

The study of lower bound of density for classical solutions can be traced back to Riemann's pioneer paper \cites{Riemann} in 1860, in which he considered a special wave interaction between two strong rarefaction waves. By studying Riemann's construction, Lipschitz continuous examples for isentropic Euler equations \eqref{p1}$\sim$\eqref{p3}
were provided in Section 82 in \cites{courant}, in which the function $\min_{x\in\mathbb{R}}\rho(x,t)$ was proved to decay to zero in an order of $O(1+t)^{-1}$ as $t\rightarrow\infty$, while the initial density is uniformly away from zero.{\footnote{The author thanks Helge Kristian Jenssen who first pointed out this result to him.}} A relative detailed discussion can be found in \cite{CPZ2},
when the adiabatic constant $\gamma=\frac{2\mathcal{N}+1}{2\mathcal{N}-1}$ with any positive integer $\mathcal{N}$.

Then there were many articles working on time-dependent lower bound on density for general classical solutions of isentropic Euler equations (\ref{p1})$\sim$(\ref{p3}) under assumption that initial density is uniformly positive. 
For rarefactive piecewise Lipschitz continuous solutions, for any $\gamma>1$,  L. Lin first proved that the density has lower bound in the order of  $O(1+t)^{-1}$ in \cite{lin2} by introducing a polygonal scheme. 
A breakthrough for general classical solutions happens in a recent paper  \cite{CPZ},
in which R. Pan, S. Zhu and the author found a lower bound of density in the order of  $O(1+t)^{-4/(3-\gamma)}$ when $1<\gamma<3$. Using this result together with Lax's decomposition in \cite{lax2}, Pan, Zhu and the author proved that gradient blowup of $u$ and/or ${\tau}$ happens in finite time if and only if the initial data are forward or backward compressive somewhere. 
Next,  for general Lipschitz continuous solution, Pan, Zhu and the author in \cite{CPZ2} improved the lower bound on density from
the order of $O(1+t)^{-4/(3-\gamma)}$  to the optimal order $O(1+t)^{-1}$
by introducing a polygonal scheme. The advantage of this method is that it works for not only classical solutions but also Lipschitz continuous solutions. And the scheme itself is of both analytical and numerical interest. However, the use of a polygonal scheme makes the proof very complex and the method seems hard to be extended to full Euler equations. 
Another result on the lower bound of density for classical solution in the order of $O(1+t)^{-1}$  when $\gamma=3$ was given by A. Bressan{\footnote{ The author knew this unpublished result through a private communication with A. Bressan.}}, where the proof relies on the study of Riccati equations established by Lax in \cite{lax2}. 

For non-isentropic full Euler equations, before this paper, the only polynomial order upper bound of $\tau$ (lower bound of $\rho$) for general classical solution was
established by Pan, Zhu and the author in \cite{CPZ}. More precisely, we showed density has a lower bound in the order of $O(1+t)^{-4/(3-\gamma)}$ when $1<\gamma<3$.

In summary, lower bound of density in optimal order $O(1+t)^{-1}$ is still not available for isentropic Euler equations with $\gamma>3$ and full nonisentropic Euler equations with $\gamma>1$, before this paper.

In this paper, we consider classical solutions of Cauchy problems of both isentropic Euler equations
and nonisentropic Euler equations. And we assume that  initial density is  uniformly positive.
We give a short proof that density 
has time-dependent lower bound in optimal order $O(1+t)^{-1}$ for  isentropic Euler equations (in Theorem \ref{thm1})  and in almost optimal order $O(1+t)^{-1-\delta}$ for any $0<\delta<\frac{1}{3}$ for full Euler equations (in Theorem \ref{thm2})  in one space dimension, respectively. 

Furthermore,  for classical solutions, we prove that $u_{x}(x,t)$  for p-system  and $\rho^\ve\, u_x$ for any $0<\ve<\frac{1}{4}$ for full Euler equations are uniformly bounded above by a constant, respectively, although they are unbounded from below when gradient blowup happens, i.e. when shock forms. In Eulerian coordinates $(y,t)$, we show  for full Euler equations
$u_{y}(y,t)$ is uniformly bounded above by a constant.

The lower bounds of density achieved in this paper can give us more precise estimate of life span of classical solution than those achieved in \cite{CPZ} and motivate us in searching lower bound of density for BV solutions including shock waves, which is a major obstacle in establishing large BV existence theory for Euler equations.
Another interesting result on a time-dependent density lower bound for isentropic Euler-Poisson equations can be found in \cite{TW} by E. Tadmor and D. Wei. 
 

The rest of the paper is divided into three sections. In Section 2, we introduce the main results and ideas in this paper. In Section 3, 
we prove Theorem \ref{thm1} for the p-system. In Section 4, we prove Theorem \ref{thm2} for the full Euler equations.

%

\section{{\bf Main results and ideas}}

Let's first introduce some variables and notations. For Euler equations (\ref{lagrangian1})$\sim$(\ref{s con}), we use variables
\beq\label{1.1}
   m\doteq e^{\frac{S}{2c_v}}\quad\hbox{and}\quad   \eta \doteq  \textstyle\frac{2\sqrt{K\gamma}}{\gamma-1}\,
\tau^{-\frac{\gamma-1}{2}}
\eeq
to take the roles of $S$ and $\tau$.
We denote the Riemann invariants 
\beq
  s \doteq u+m\,\eta\quad\hbox{and}\quad r \doteq u-m\,\eta\,
\label{1.2}
\eeq
respectively,
and gradient variables
\beq\label{1.3}
\alpha \doteq u_x+m\,\eta_x+ \textstyle\frac{\gamma-1}{\gamma}m_x\, \eta\quad\hbox{and}\quad
\beta \doteq u_x-m\,\eta_x- \textstyle\frac{\gamma-1}{\gamma}m_x\, \eta\,.
\eeq

For the isentropic Euler equations (p-system) (\ref{p1})$\sim$(\ref{p3}), whose solutions are special solutions of full Euler equations \eqref{lagrangian1}$\sim$\eqref{introduction 3} when we restrict our consideration on the classical solution, Riemann invariants are
\beq\label{srdef}
s= u+\eta \quad\hbox{and}\quad  r= u-\eta
\eeq
and
\beq\label{alpha} \left.\begin{array}{l}
\alpha=u_x+\eta_x=s_x\quad\hbox{and}\quad\beta=u_x-\eta_x= r_x\,.\end{array}\right.\eeq 

The main results in this paper are listed in 
the following two theorems: Theorem \ref{thm1} for p-system and Theorem \ref{thm2} for full Euler equations.
\begin{theorem}\label{thm1}
Let $\big(u(x,t), \tau(x,t)\big)$ be a  $C^1$ solution of the isentropic Euler equations (\ref{p1})$\sim$(\ref{p3})
in the region $(x,t)\in \mathbb R\times[0,T)$, where $T$ can be any finite positive constant or infinity.
Assume
$u(x,0)$, $\tau(x,0)>0$, $\rho(x,0)=1/\tau(x,0)$, $\alpha(x,0)$ and $\beta(x,0)$ are all uniformly bounded,
where $\alpha$ and $\beta$ take the form in \eqref{alpha}.

Let $M$ be an upper bound of $\alpha(x,0)$ and $\beta(x,0)$, i.e.
\beq\label{al-be-0}
\max_{x\in\mathbb R}\Big\{\alpha(x,0),\beta(x,0)\Big\}<M
\eeq
then 
\beq\label{al-be-t}
\max_{(x,t)\in\mathbb R\times[0,T)}\Big\{\alpha(x,t),\beta(x,t)\Big\}<M\,.
\eeq
This gives
\beq\label{uxtautM}
\max_{(x,t)\in\mathbb R\times[0,T)}\big\{\tau_t\big\}=\max_{(x,t)\in\mathbb R\times[0,T)}\big\{u_x\big\}<M
\eeq
by (\ref{alpha}) and (\ref{p1}).
Hence, there exist positive constants $M_1$ and $M_2$ independent of $T$ such that
\beq\label{main1}
\min_{x}{\rho(x,t)}\geq\frac{M_1}{M_2+t}\,.
\eeq

\end{theorem}

The key step in the proof of Theorem \ref{thm1} is to prove (\ref{al-be-t}).
In fact, suppose (\ref{al-be-t}) is correct, then by the conservation of mass (\ref{p1}) and (\ref{alpha}),
we can easily prove \eqref{uxtautM}:
\beq\label{main0}
\tau_t=u_x=\frac{1}{2}(\alpha+\beta)<M
\eeq
which directly gives (\ref{main1}), together with the initial condition. To prove (\ref{al-be-t}), we need to study the characteristic decomposition established by Lax
in \cite{lax2}. 
The key idea is to find an invariant domain on $\alpha$ and $\beta$.

One conclusion that we can draw from (\ref{al-be-0})$\sim$(\ref{al-be-t}) is that although
the variables $\alpha$ and  $\beta$ might increase along forward and backward characteristics, respectively, the function $\max_{x\in \mathbb R}\{\alpha(x,t),\beta(x, t)\}$ is not increasing
with respect to $t$, which means that the maximum rarefaction of classical solution is not increasing. This result can be easily seen from the fact that 
(\ref{al-be-t}) is still correct if we change $0$ in 
(\ref{al-be-0}) into any $t^*\in(0,t)$.

\begin{remark}
Under assumptions in Theorem \ref{thm1},
in Eulerian coordinates $(y, t),$ the inequality (\ref{uxtautM}) gives that smooth solutions in 
the region $(y,t)\in\mathbb R\times[0,T)$ satisfy
\beq\label{remark1a}
\max_{(y,t)\in\mathbb R\times[0,T)}\Big\{\frac{ u_{y}}{\rho}\Big\}<M\,,
\eeq
where $M$ is the constant given in \eqref{al-be-0},
because $\rho\, u_x(x,t)=u_{y}(y,t)$. See \cite{smoller} for the transformation between Eulerian and Lagrangian coordinates.

Since $\rho$ is uniformly bounded above, which can be easily proved by the fact that Riemann invariants $s$ and $r$ are initially bounded and are constant along forward and backward characteristics, respectively,
we know
\[
\max_{(y,t)\in\mathbb R\times[0,T)}\Big\{ u_{y}\Big\}<\bar{M}\,,
\] 
for some constant $\bar{M}$ independent of $T$. 

\end{remark}

\bigskip
Then we consider the full Euler equations.
\begin{theorem}\label{thm2}
Let $\big(u(x,t), \tau(x,t), S(x)\big)$ be a $C^1$ solution 
of full Euler equations \eqref{lagrangian1}$\sim$\eqref{introduction 3}
in the region $(x,t)\in \mathbb R\times[0,T)$. Here, $T$ can be any finite positive constant or infinity. Assume that initial data
$u(x,0)$, $\tau(x,0)>0$, $\rho(x,0)=1/\tau(x,0)$, $S(x)$, $S'(x)$, $\alpha(x,0)$ and $\beta(x,0)$ are all uniformly bounded and total variation of $S(x)$ is finite, where $\alpha$ and $\beta$ satisfy (\ref{1.3}).
Then, for any 
\[0<\ve<\frac{1}{4}\,,\] there exists constant $N_0$ independent of $T$, such that 
\beq\label{uxtautM2}
\max_{(x,t)\in\mathbb R\times[0,T)}\Big\{\rho^\ve\cdot\tau_t\Big\}=\max_{(x,t)\in\mathbb R\times[0,T)}\Big\{\rho^\ve\cdot u_x\Big\}<N_0\,,
\eeq
and there exist positive constants $N_1$ and $N_2$ independent of $T$ such that
\beq\label{main2}
\min_{x}{\rho(x,t)}\geq\Big(\frac{N_1}{N_2+t}\Big)^{1+\delta}\,,
\eeq
where $\delta= \frac{\ve}{1-\ve}>0$.

\end{theorem}

We first prove a result in Lemma \ref{lemma_4.1} taking the similar role as (\ref{al-be-t}) in Theorem \ref{main1}. 
In fact, we find uniform bounds on gradient variables $\rho^{\ve}\alpha$ and $\rho^{\ve}\beta$, using which we can easily prove (\ref{uxtautM2}) by \eqref{1.3} and \eqref{lagrangian1}:
\[
\rho^{\ve}\tau_t=\rho^{\ve}u_x=\frac{1}{2}(\rho^{\ve}\alpha+\rho^{\ve}\beta)<\hbox{Constant}\,,
\] 
then show \eqref{main2}.
The reason why we use $\rho^{\ve}\alpha$ and $\rho^{\ve}\beta$  instead of $\alpha$ and $\beta$ is to control the lower order terms in the Riccati equations produced by the varying entropy. 
The proof of Theorem \ref{thm2} also relies on the uniform constant upper bound of density established in \cite{G8} by R. Young, Q. Zhang and the author for classical solutions when total variation of initial entropy
is finite.

\begin{remark}
Under assumptions in Theorem \ref{thm2},
in Eulerian coordinates $(y, t),$ the inequality (\ref{uxtautM}) gives that the classical solution in 
the region $(y,t)\in\mathbb R\times[0,T)$ satisfies
\[
\max_{(y,t)\in\mathbb R\times[0,T)}\Big\{\frac{ u_{y}}{\rho^{1-\ve}}\Big\}<N_0\,.
\]

Since $\rho$ is uniformly bounded above under assumptions in Theorem \ref{thm2},
we know
\[
\max_{(y,t)\in\mathbb R\times[0,T)}\Big\{ u_{y}\Big\}<\bar{N_0}\,,
\] 
for some constant $\bar{N_0}$ independent of $T$.

See \cite{smoller} for the transformation between Eulerian and Lagrangian coordinates.
Since this result is a local result, we only need to assume that initial entropy is locally BV.
\end{remark}

One direct application of Theorem \ref{thm2} is that one can use (\ref{main2}) to improve the life-span estimate established in \cite{CPZ} when $1<\gamma<3$ which depends on the time-dependent lower bound of density. We leave this to the reader.

 %
\section{{\bf Lower bound of density for p-system: The proof of Theorem \ref{thm1}}}
We first introduce the characteristic decompositions for $C^1$ solution of p-system.
For any classical solution for (\ref{p1})$\sim$(\ref{p3}), the Riemann invariants $s$ and $r$ in (\ref{srdef}) are constant along forward and backward characteristics, respectively, 
\beq\label{srcon}
 \partial_+s=0\quad\hbox{and}\quad  \partial_- r=0
\eeq
with
\[
 \partial_+=\partial_t+c\partial_x\quad\hbox{and}\quad  \partial_-=\partial_t-c\partial_x\,
\]
and wave speed
\[
  c=\sqrt{-p_\tau}=
  \sqrt{K\,\gamma}\,{\tau}^{-\frac{\gamma+1}{2}}\,.
\]
Furthermore, gradient variables $\alpha=s_x$ and $\beta=r_x$ defined in \eqref{alpha} satisfy the following Riccati equations.
{\proposition\label{remark}\cite{G3} The classical solution
in (\ref{p1})$\sim$(\ref{p3}) satisfy
 
 \beq
\left.\begin{array}{l}\partial_+\alpha=k_1\{\alpha\beta-\alpha^2\}\,, \label{rem1}
\end{array}\right.\eeq
and \beq \left.\begin{array}{l}\partial_-\beta=k_1\{\alpha\beta-\beta^2\}\,,\label{rem2}
\end{array}\right.\eeq
where \beq\left.\begin{array}{l}
k_1\doteq \frac{(\gamma+1)K_c}{2(\gamma-1)}  \eta^{\frac{2}{\gamma-1}}, \label{k1 def}\end{array}\right.\eeq
where $K_c$ is a positive constant given in \eqref{Kdefs}. The function $\eta>0$ is defined in \eqref{1.1}.}

Equations \eqref{rem1} and \eqref{rem2} are special examples of Lax's decompositions
 in \cite{lax2} 
for general hyperbolic systems with two unknowns.
See detailed derivation of \eqref{rem1} and \eqref{rem2} in \cite{G3}.

\begin{figure}[htp] \centering
		\includegraphics[width=.25\textwidth]{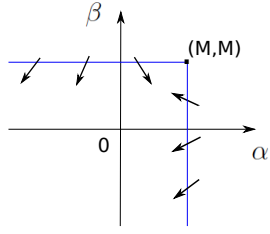}
		\caption{$\max\{\alpha,\beta\}<M$ is an invariant domain.
		Note: $\alpha$ (or $\beta$) might increase along characteristic.\label{key}}
	\end{figure}

\begin{remark}
The idea for the proof of \eqref{al-be-t} can be seen from Figure \ref{key}.

\end{remark}

Before the proof, we remark on one fact that $\rho$, $\eta$, $c$ and $k_1$ are all bounded above
by some constants if assumptions in Theorem \ref{thm1} are satisfied.
This can be easily obtained by (\ref{srcon}), which says that $s$ and $r$ are constant along forward and backward characteristics. As a consequence, $\rho$, $\eta$, $c$ and function $k_1$ are all uniformly bounded from above.
Denote 
\beq\label{kK1}
K_1\doteq\max_{(x,t)\in\mathbb R\times [0,T)}k_1(x,t)\,,
\eeq
where $K_1$ is a constant only depending on $\gamma$ and initial condition.

\bigskip

\begin{figure}[htp] \centering
		\includegraphics[width=.4\textwidth]{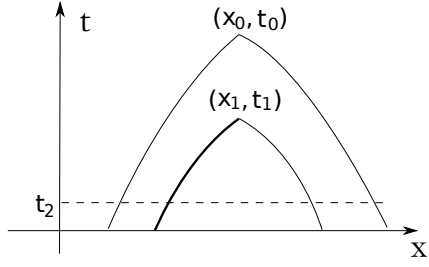}
		\caption{Proof of Theorems \ref{thm1} and  \ref{thm2}.\label{fig0}}
	\end{figure}

\paragraph{{\bf Proof of Theorem \ref{thm1}}}
We first prove  (\ref{al-be-t}) by contradiction. Without loss of generality, assume that 
$\alpha(x_0,t_0)=M$ at some point $(x_0,t_0)$. See Figure \ref{fig0}.

Because wave speed $c$ is bounded above, we can find the characteristic triangle with vertex $(x_0,t_0)$ and lower boundary
on the initial line $t=0$, denoted by $\Omega$. 

Then we can find the first time $t_1$ such that $\alpha=M$ or $\beta=M$ in $\Omega$. More precisely,
\[
\max_{(x,t)\in \Omega,\, t<t_1}\Big(\alpha(x,t),\beta(x,t)\Big)<M,\]
and $\alpha(x_1,t_1)=M$ or/and $\beta(x_1,t_1)=M$ for some $(x_1,t_1)\in\Omega$. Without loss of generality, still assume 
$\alpha(x_1,t_1)=M$. The proof for another case is entirely same.
Let's denote the characteristic triangle with vertex $(x_1,t_1)$ as  $\Omega_1\in\Omega$, then
\beq\label{p_main_1}
\max_{(x,t)\in \Omega_1,\, t<t_1}\Big(\alpha(x,t),\beta(x,t)\Big)<M,\eeq
and $\alpha(x_1,t_1)=M$.
By the continuity of $\alpha$,  we could find a time $t_2\in[0,t_1)$ such that,
\beq\label{p_main_2}
\alpha(x,t)>0\,,\quad \hbox{for any}\quad (x,t)\in \Omega_1\quad \hbox{and}\quad t\geq t_2\,.
\eeq

Next we derive a contradiction. By \eqref{rem1}, \eqref{kK1} and 
\eqref{p_main_1}$\sim$\eqref{p_main_2}, along the forward characteristic segment  through $(x_1,t_1)$ when $t_2\leq t<t_1$,
\begin{eqnarray*}
\partial_+\alpha&=&k_1\{\alpha\beta-\alpha^2\}\leq K_1\{M\alpha-\alpha^2\}
\end{eqnarray*}
which gives, through integration along characteristic,
\begin{eqnarray*}
&\frac{d\alpha}{(M-\alpha)\alpha}&\leq\quad \textstyle K_1dt\\
\Rightarrow&\frac{1}{M}\ln \frac{\alpha(t)}{M-\alpha(t)}&\leq\textstyle\quad\frac{1}{M}\ln \frac{\alpha(t_2)}{M-\alpha(t_2)}+K_1 (t-t_2)\,.
\end{eqnarray*}
As $t\rightarrow t_1-$, left hand side approaches infinity 
while right hand side approaches a finite number, which gives a contradiction. Hence we prove that (\ref{al-be-t}) is correct, i.e. $\alpha$ and $\beta$
are uniformly bounded above. Then by the conservation of mass (\ref{p1}) and (\ref{alpha}), we have
\eqref{main0} then \eqref{uxtautM}, which directly gives \eqref{main1}.
Hence we complete the proof of Theorem \ref{thm1}.

\begin{remark} \label{remark_gp}
Theorem \ref{thm1} can be extended to the case with general pressure law $p=p(\tau)$ with $p_\tau<0$, $p_{\tau\tau}>0$ and some other suitable
conditions on $p$. We leave this to the reader and refer the reader  
to \cite{G5} for the Riccati equations and the definitions of $\alpha$ and $\beta$.
For full Euler equations, the extension of Theorem \ref{thm2} to general pressure law is still not available because the current result on uniform upper bound of density is only available for $\gamma$-law pressure.
\end{remark}

\section{{\bf Full compressible Euler equations}}
\subsection{Equations and coordinates}
We first introduce some notations and existing equations for $C^1$ solutions
of full Euler equations \eqref{lagrangian1}$\sim$\eqref{introduction 3}.
Recall we use new variables $m$ and $\eta$ to take the roles of $S$ and $\tau$, respectively:
\beq
   m=e^{\frac{S}{2c_v}}\label{m def}
\eeq
and 
\beq
   \eta  = \TS\frac{2\sqrt{K\gamma}}{\gamma-1}\,
\tau^{-\frac{\gamma-1}{2}}\,.\label{z def}
\eeq
Without confusion, we still use $c$ to denote  the nonlinear Lagrangian wave speed for full Euler equations, where
\beq
  c=\sqrt{-p_\tau}=
  \sqrt{K\,\gamma}\,{\tau}^{-\frac{\gamma+1}{2}}\,e^{\frac{S}{2c_v}}\,.
\label{c_def2}
\eeq
The forward and backward characteristics are described by
\beq\label{pmc_full}
  \frac{dx}{dt}=c \com{and} \frac{dx}{dt}=-c\,,
\eeq
and we denote the corresponding directional derivatives along these characteristics by
\[
  \pp := \dbyd t+c\,\dbyd x \com{and}
  \pn := \dbyd t-c\,\dbyd x\,,
\]
respectively.

It follows that
\begin{align}
  \tau&=K_{\tau}\,\eta^{-\frac{2}{\gamma-1}}\,,\nn\\
  p&=K_p\, m^2\, \eta^{\frac{2\gamma}{\gamma-1}}\,,\label{tau p c}\\
  c&=c(\eta,m)=K_c\, m\, \eta^{\frac{\gamma+1}{\gamma-1}}\,.\nn
\end{align}
with positive constants
\beq
  K_\tau:=\Big(\frac{2\sqrt{K\gamma}}{\gamma-1}\Big)^\frac{2}{\gamma-1}\,,
\quad
  K_p:=K\,K_\tau^{-\gamma},\com{and}
  K_c:=\TS\sqrt{K\gamma}\,K_\tau^{-\frac{\gamma+1}{2}}\,,
\label{Kdefs}
\eeq
so that also
\beq
   K_p=\TS\frac{\gamma-1}{2\gamma}K_c \com{and}
   K_\tau K_c=\frac{\gamma-1}{2}\,.
\label{KpKcRela}
\eeq

In these coordinates, for $C^1$ solutions, equations
\eq{lagrangian1}--\eqref{introduction 3} are equivalent to
\begin{align}
  \eta_t+\frac{c}{m}\,u_x&=0\,, \label{lagrangian1 zm}\\
  u_t+m\,c\,\eta_x+2\frac{p}{m}\,m_x&=0\,,\label{lagrangian2 zm}\\
  m_t&=0\label{lagrangian3 zm}\,,
\end{align}
where the last equation comes from (\ref{s con}), which is equivalent to
(\ref{lagrangian3}), c.f. \cite{smoller}.  Note that, while the solution remains $C^1$,
$m=m(x)$ is given by the initial data and can be regarded as a
stationary quantity.

Recall that we denote the Riemann invariants by
\beq
  r:=u-m\,\eta\quad\hbox{and}\quad s:=u+m\,\eta\,.
\label{r_s_def}
\eeq
Different from the isentropic case ($m$ constant), for general non-isentropic flow, $s$ and $r$
vary along characteristics.  
Also recall we denote gradient variables
\beq \left.\begin{array}{l}
\alpha=u_x+m\eta_x+\frac{\gamma-1}{\gamma}m_x
\eta,\end{array}\right.\label{al_def}\eeq \beq
\left.\begin{array}{l}\beta=u_x-m\eta_x-
\frac{\gamma-1}{\gamma}m_x \eta,\end{array}\right.\label{be_def} \eeq
which satisfy the following Riccati equations. See detailed derivation in  \cite{G3}.
{\begin{proposition}\cite{G3}\label{remark2} The classical solutions
for (\ref{lagrangian1})$\sim$(\ref{lagrangian3}) satisfy \beq
\left.\begin{array}{l}\partial_+\alpha=k_1\{k_2
(3\alpha+\beta)+\alpha\beta-\alpha^2\}, \label{frem1}
\end{array}\right.\eeq
and \beq \left.\begin{array}{l}\partial_-\beta=k_1\{-k_2
(\alpha+3\beta)+\alpha\beta-\beta^2\},\label{frem2}
\end{array}\right.\eeq
where \beq\left.\begin{array}{l}
k_1=\frac{(\gamma+1)K_c}{2(\gamma-1)} \eta^{\frac{2}{\gamma-1}}, \quad
k_2=\frac{\gamma-1}{\gamma(\gamma+1)}\eta\, m_x. \label{k def}\end{array}\right.\eeq
\end{proposition}
Proposition \ref{remark} is in fact a corollary of Proposition \ref{remark2} for the isentropic case
in which $m_x\equiv0$.

%
\subsection{Uniform upper bound on density}
In this part, we review a result on the uniform upper bounds of $|u|$ and $\rho$  established by the author, R. Young and Q. Zhang in \cite{G6}, for later references.

In this section, we always assume all initial conditions in Theorem \ref{thm2} are satisfied.
So that
\beq
  V := \frac{1}{2c_v}\int_{-\infty}^{+\infty}|S'(x)|\;dx
     = \int_{-\infty}^{+\infty}\frac{|m'(x)|}{m(x)}\;dx<\infty\,,
\label{Vdef}
\eeq
while also, by \eq{m def},
\beq
  0 < M_L < m(\cdot) < M_U\,,
\label{m_bounds}
\eeq
for some constants $M_L$ and $M_U$.  Also there exist positive constants $M_s$ and $M_r$,
such that, in the initial data,
\beq
  |s(\cdot,0)|<M_s \com{and}
  |r(\cdot,0)|<M_r\,.
\label{Mrs}
\eeq
 
In the following proposition established in \cite{G6}, $|u|$ and $\rho$ are shown to be uniformly bounded above.

\begin{proposition}{\em \cite{G6}}
\label{Thm_upper}Assume all initial conditions in Theorem \ref{thm2} are satisfied.
And assume system \eqref{lagrangian1}$\sim$\eqref{introduction 3}  has a  $C^1$ solution when $t\in[0,T)$, then one
has the uniform bounds
\beq\label{u_rho_bounds}
   |u(x,t)|\leq\frac{L_1+L_2}{2}{M_U}^{\frac{1}{2\gamma}}
\com{and}
   \eta(x,t)\leq\frac{L_1+L_2}{2}{M_L}^{\frac{1}{2\gamma}-1},
\eeq
where
\begin{align*}
  L_1 &:= M_s+\ol V\,M_r+\ol V\,(\ol V\,M_s+{\ol V}^2\,M_r)
	\,e^{{\ol V}^2},\\
  L_2 &:= M_r+\ol V\,M_s+\ol V\,(\ol V\,M_r+{\ol V}^2\,M_s)
	\,e^{{\ol V}^2},
\end{align*}
and
\[\ol V := \frac{V}{2\gamma}\,.
\] 
Constants $L_1$ and $L_2$ both clearly depend only on the
initial data and $\gamma$.
Here $T$ can be any positive number or infinity. And the bounds are independent of $T$.
\end{proposition}

%
\subsection{Proof of Theorem \ref{thm2}}
Similar as Theorem \ref{thm1} for p-system, the key idea is still to get the uniform upper bound of some gradient variables measuring rarefaction.

However, we cannot directly get the  uniform upper bound of $\alpha$ and $\beta$. In fact,
comparing to (\ref{rem1})$\sim$(\ref{rem2}),  equations  (\ref{frem1})$\sim$(\ref{frem2})
include some first order terms in the right hand side.  In order to cope with them, we introduce some
new gradient variables
\beq\label{albeve}
\alpha_\ve=\eta^{\frac{2\ve}{\gamma-1}}\, \alpha\quad\hbox{and}
\quad
\beta_\ve=\eta^{\frac{2\ve}{\gamma-1}}\, \beta.
\eeq
Using (\ref{lagrangian1 zm}), we have
\[
 \partial_+ \eta=\eta_t+c\eta_x=-\frac{c}{m}u_x+c\eta_x=-K_c\eta^{\frac{\gamma+1}{\gamma-1}}\beta
 -\textstyle\frac{\gamma-1}{\gamma}K_c \eta^{\frac{2\gamma}{\gamma-1}}m_x \,,
\]
and
\[
 \partial_- \eta=\eta_t-c\eta_x=-\frac{c}{m}u_x-c\eta_x=-K_c\eta^{\frac{\gamma+1}{\gamma-1}}\alpha
 +\textstyle\frac{\gamma-1}{\gamma}K_c \eta^{\frac{2\gamma}{\gamma-1}}m_x \,,
\]
then it is easy to prove the next lemma by Proposition \ref{remark2}.
{\begin{lemma}\label{remark3} The classical solutions
in (\ref{lagrangian1})$\sim$(\ref{lagrangian3}) satisfy
\beq
\left.\begin{array}{l}\partial_+\alpha_\ve=k_{1\ve}\left\{k_{2\ve}
(3\alpha_\ve-4\ve\alpha_\ve+\beta_\ve)+(1-\frac{4\ve}{\gamma+1})\alpha_\ve\beta_\ve-\alpha_\ve^2\right\}, \label{rem3}
\end{array}\right.\eeq
and \beq \left.\begin{array}{l}\partial_-\beta_\ve=k_{1\ve}\{-k_{2\ve}
(\alpha_\ve+3\beta_{\ve}-4\ve\beta_\ve)+(1-\frac{4\ve}{\gamma+1})\alpha_\ve\beta_\ve-\beta_\ve^2\},\label{rem4}
\end{array}\right.\eeq
where \beq\left.\begin{array}{l}
k_{1\ve}=\frac{(\gamma+1)K_c}{2(\gamma-1)} \eta^{\frac{2}{\gamma-1}(1-\ve)}, \quad
k_{2\ve}=\frac{\gamma-1}{\gamma(\gamma+1)}\eta^{1+\frac{2}{\gamma-1}\ve}\, m_x, \label{fk def}\end{array}\right.\eeq
and 
\beq\label{vecon}
0<\ve<\frac{1}{4}.
\eeq
\end{lemma}}

Note, for any $C^1$ solutions in $(x,t)\in \mathbb R\times[0,T)$ satisfying initial conditions in 
Theorem \ref{thm2}, using Proposition \ref{Thm_upper},  for any $\ve$ satisfying \eqref{vecon},
we know $|k_{1\ve}(x,t)|$ and
$|k_{2\ve}(x,t)|$ are both uniformly bounded above:
\beq\label{k_ve}
|k_{1\ve}(x,t)|< \hat K_{1}\quad \hbox{and}\quad |k_{2\ve}(x,t)|< \hat K_{2},
\eeq
where constants $\hat K_{1}$ and $\hat K_{2}$ only depend on initial conditions and $\gamma$
but independent of $\ve$.

Next we give the key lemma which will be proved later.

{\begin{lemma}\label{lemma_4.1}
Suppose the initial conditions in Theorem \ref{thm2} are satisfied.  For any $\ve$ satisfying
\eqref{vecon}, let $N$ be an upper bound of $\alpha_\ve(x,0)$ and $\beta_\ve(x,0)$, i.e.
\beq\label{2al-be-0}
\max_{x\in\mathbb R}\Big\{\alpha_\ve(x,0),\beta_\ve(x,0)\Big\}<N
\eeq
where constant $N$ also satisfies
\beq\label{M_con}
N>\max\big\{\textstyle\frac{4(\gamma+1)\hat K_{2}}{\ve},\frac{2\hat K_{2}}{1-\frac{4\ve}{\gamma+1}} \big\},
\eeq
then 
\beq\label{2al-be-t}
\max_{(x,t)\in\mathbb R\times[0,T)}\Big\{\alpha_\ve(x,t),\beta_\ve(x,t)\Big\}<N\,.
\eeq
\eqref{vecon}.
\end{lemma}}

%
%

\paragraph{{\bf Proof of Theorem \ref{thm2}}}
We only have to show Lemma \ref{lemma_4.1}. In fact, if  
Lemma \ref{lemma_4.1} is proved, then by the conservation of mass \eqref{lagrangian1} and definitions of $\alpha_\ve$ and $\beta_\ve$ in \eqref{albeve} and \eqref{al_def}$\sim$\eqref{be_def}, we have
\[
\eta^{\frac{2\ve}{\gamma-1}}\tau_t=\eta^{\frac{2\ve}{\gamma-1}} u_x
=\frac{1}{2}(\alpha_\ve+\beta_\ve)<N
\]
which gives that, by \eqref{z def}, $\tau=1/\rho$ and initial density has positive lower bound,
there exists positive constants $N_1$ and $N_2$, such that
\[
\rho>(\frac{N_1}{N_2+t})^{1+\delta}
\]
where
\[
\delta=\frac{\ve}{1-\ve}\,.
\]
Then it is easy to see that all results in
Theorem \ref{thm2} are correct.

Now we prove  Lemma \ref{lemma_4.1} by contradiction. We still use Figure \ref{fig0}. Without loss of generality, assume that 
$\alpha_\ve(x_0,t_0)=N$, at some point $(x_0,t_0)$.

Because wave speed $c$ is bounded above, we can find the characteristic triangle with vertex $(x_0,t_0)$ and lower boundary
on the initial line $t=0$, denoted by $\Omega$. 

Then we can find the first time $t_1$ such that $\alpha_\ve=N$ or $\beta_\ve=N$ in $\Omega$. More precisely,
\[
\max_{(x,t)\in \Omega,\, t<t_1}\Big(\alpha_\ve(x,t),\beta_\ve(x,t)\Big)<N,\]
and $\alpha_\ve(x_1,t_1)=N$ or/and $\beta_\ve(x_1,t_1)=N$ for some $(x_1,t_1)\in\Omega$. Without loss of generality, still assume 
$\alpha_\ve(x_1,t_1)=N$. The proof for another case is entirely same.
Let's denote the characteristic triangle with vertex $(x_1,t_1)$ as  $\Omega_1\in\Omega$, then
\[
\max_{(x,t)\in \Omega_1,\, t<t_1}\Big(\alpha_\ve(x,t),\beta_\ve(x,t)\Big)<N,\]
and $\alpha_\ve(x_1,t_1)=N$.

Then we divide the problem into two cases:
\begin{itemize}
\item[I.] $N\geq\beta_\ve(x_1,t_1)>-\frac{N}{2}$;
\item[II.] $\beta_\ve(x_1,t_1)\leq-\frac{N}{2}$.
\end{itemize}

In case I, by the continuity of $\alpha_\ve$ and $\beta_\ve$ and our construction,  we can find a time $t_2\in[0,t_1)$ such that,
\beq\label{final1}
\frac{N}{2}<\alpha_\ve(x,t)<N\quad\hbox{and}\quad |\beta_\ve|< N\,,\quad \hbox{for any}\quad (x,t)\in \Omega_1\quad \hbox{and}\quad t_2\leq t< t_1\,.
\eeq
Then using \eqref{rem3}, \eqref{k_ve}, \eqref{M_con} and (\ref{final1}), along the forward characteristic segment through $(x_1,t_1)$, when $t_2\leq t<t_1$, we have
\begin{eqnarray*}
\partial_+\alpha_\ve\leq\textstyle k_{1\ve}(1-\frac{4\ve}{\gamma+1})\,(\alpha_\ve\beta_\ve-\alpha_\ve^2)\leq\tilde{K}_{1}\,(N\alpha_\ve-\alpha_\ve^2)
\end{eqnarray*}
with 
\[\tilde{K}_{1}\doteq \hat K_{1}(1-\frac{4\ve}{\gamma+1}),\]
which gives, through integration along characteristic,
\begin{eqnarray*}
&\frac{d\alpha_\ve}{(N-\alpha_\ve)\alpha_\ve}&\textstyle\leq\quad\tilde{K}_{1} dt\\
\Rightarrow&\frac{1}{N}\ln \frac{\alpha_\ve(t)}{N-\alpha_\ve(t)}&\leq\textstyle\quad\frac{1}{N}\ln \frac{\alpha_\ve(t_2)}{N-\alpha_\ve(t_2)}+\tilde{K}_{1} (t-t_2)\,.
\end{eqnarray*}
As $t\rightarrow t_1-$, left hand side approaches infinity 
while right hand side approaches a finite number, which gives a contradiction.

In case II, by the continuity of $\alpha_\ve$,  we could find a time $t_3\in[0,t_1)$ such that,
\beq\label{final2}
\frac{N}{2}<\alpha_\ve(x,t)<N\quad\hbox{and}\quad \beta_\ve(x,t)< -\frac{N}{4}\,,\quad \hbox{for any}\quad (x,t)\in \Omega_1\quad \hbox{and}\quad t_3\leq t< t_1\,.
\eeq
which gives, by \eqref{M_con},
\[
\textstyle\big(k_{2\ve}+(1-\frac{4\ve}{\gamma+1})\alpha_\ve\big)\,\beta_\ve<0\,.
\]
Hence by \eqref{M_con}, \eqref{vecon} and \eqref{final2}, we have 
\[\partial_+\alpha_\ve<k_{1\ve}\left\{k_{2\ve}
(3-4\ve)\alpha_\ve-\alpha_\ve^2\right\}<0\,.
\]
As a consequence, $\alpha_\ve$ decreases on $t$ along the forward characteristic line through $(x_1,t_1)$,
when $t_3\leq t<t_1$,
which contradicts to that $\alpha_\ve(x_1,t_1)=N$ while $\alpha_\ve(x,t)<N$
when $(x,t)\in \Omega_1$ and $t_3\leq t< t_1$. Hence Lemma \ref{lemma_4.1} is proved.
This completes the proof of Theorem \ref{thm2}.

%

%


\end{document}